\documentclass[12pt]{article}
\usepackage{amsmath}
\usepackage{amssymb}
\usepackage{amsfonts}
\usepackage{nccmath}
\usepackage{comment}
\usepackage{color}
\usepackage{mathtools}
\usepackage{tikz}
\usepackage{enumitem}

\DeclarePairedDelimiter{\floor}{\lfloor}{\rfloor}
\usepackage{amsthm}

\newtheorem*{remark*}{Remark}
\newtheorem*{notation*}{Notation}
\newtheorem{theorem}{Theorem}[section]
\newtheorem{definition}[theorem]{Definition}
\newtheorem{proposition}[theorem]{Proposition}
\newtheorem{problem}[theorem]{Problem}
\newtheorem{lemma}[theorem]{Lemma}
\newtheorem{corollary}[theorem]{Corollary}
\newtheorem{conjecture}[theorem]{Conjecture}

\numberwithin{equation}{section}
\DeclareMathOperator{\sat}{sat}
\DeclareMathOperator{\ssat}{ssat}
\DeclareMathOperator{\ex}{ex}

\def \F {\mathcal{F}}

\tikzstyle{cir} = [draw, circle, minimum height= 20 mm]

\title{Minimizing the numbers of cliques and cycles of fixed size in an $F$-saturated graph}
\author{Debsoumya Chakraborti\thanks{Department of Mathematical Sciences, Carnegie Mellon University. Email: {\tt dchakrab@cmu.edu}. Research supported in part by National Science Foundation CAREER Grant DMS-1455125.} \ and Po-Shen Loh\thanks{Department of Mathematical Sciences, Carnegie Mellon University. Email: {\tt ploh@cmu.edu}. Research supported in part by National Science Foundation CAREER Grant DMS-1455125.}}
\begin{document}
\maketitle
\begin{abstract}
This paper considers two important questions in the well-studied theory of graphs that are $F$-saturated. A graph $G$ is called $F$-saturated if $G$ does not contain a subgraph isomorphic to $F$, but the addition of any edge creates a copy of $F$. We first resolve a fundamental question of minimizing the number of cliques of size $r$ in a $K_s$-saturated graph for all sufficiently large numbers of vertices, confirming a conjecture of Kritschgau, Methuku, Tait, and Timmons. We also go further and prove a corresponding stability result. Next we minimize the number of cycles of length $r$ in a $K_s$-saturated graph for all sufficiently large numbers of vertices, and classify the extremal graphs for most values of $r$, answering another question of Kritschgau, Methuku, Tait, and Timmons for most $r$. 

We then move on to a central and longstanding conjecture in graph saturation made by Tuza, which states that for every graph $F$, the limit $\lim_{n \rightarrow \infty} \frac{\sat(n, F)}{n}$ exists, where $\sat(n, F)$ denotes the minimum number of edges in an $n$-vertex $F$-saturated graph. Pikhurko made progress in the negative direction by considering families of graphs instead of a single graph, and proved that there exists a graph family $\mathcal{F}$ of size $4$ for which $\lim_{n \rightarrow \infty} \frac{\sat(n, \mathcal{F})}{n}$ does not exist (for a family of graphs $\mathcal{F}$, a graph $G$ is called $\mathcal{F}$-saturated if $G$ does not contain a copy of any graph in $\mathcal{F}$, but the addition of any edge creates a copy of a graph in $\mathcal{F}$, and $\sat(n, \mathcal{F})$ is defined similarly). We make the first improvement in 15 years by showing that there exist infinitely many graph families of size $3$ where this limit does not exist. Our construction also extends to the generalized saturation problem when we minimize the number of fixed-size cliques. We also show an example of a graph $F_r$ for which there is irregular behavior in the minimum number of $C_r$'s in an $n$-vertex $F_r$-saturated graph.
\end{abstract}

\section{Introduction}
Extremal graph theory focuses on finding the extremal values of certain parameters of graphs under certain natural conditions. One of the most well-studied conditions is $F$-freeness. For graphs $G$ and $F$, we say that $G$ is $F$-free if $G$ does not contain a subgraph isomorphic to $F$. This gives rise to the most fundamental question of finding the Tur\'an number $\ex(n, F)$, which asks for the maximum number of edges in an $n$-vertex $F$-free graph. The asymptotic answer is known for most graphs $F$, with the exception of bipartite $F$ where the most intricate and unsolved cases appear (see, e.g., \cite{FS} and \cite{S} for nice surveys). Recently, Alon and Shikhelman \cite{AS} introduced a natural generalization of the Tur\'an number. They systematically studied $\ex(n, H, F)$, which denotes the maximum number of copies of $H$ in an $n$-vertex $F$-free graph. Note that the case $H = K_2$ is the standard Tur\'an problem, i.e., $\ex(n, K_2, F) = \ex(n, F)$.

While the Tur\'an number asks for the maximum number of edges in an $F$-free graph, another very classical problem concerns the minimum number of edges in an $F$-free graph with a fixed number of vertices. This problem is not interesting as stated because the empty graph is the obvious answer. In much of the research, this issue is resolved by imposing the additional condition that adding any edge to $G$ will create a copy of $F$. With this additional condition, we say that $G$ is $F$-saturated. A moment's thought will convince the reader that when maximizing the number of edges, this additional condition does not change the problem at all. On the other hand, this new condition makes the edge minimization problem very interesting, and this area of research is commonly known as graph saturation. Let the saturation function $\sat(n, F)$ denote the minimum number of edges in an $n$-vertex $F$-saturated graph. Erd\H{o}s, Hajnal, and Moon \cite{EHM} started the investigation of this area with the following beautiful result.

\begin{theorem} [Erd\H{o}s, Hajnal, and Moon 1964] \label{EHM}
For every $n \ge s \ge 2$, the saturation number $$\sat(n, K_s) = (s -2)(n-s+2) + \binom{s-2}{2}.$$ Furthermore, there is a unique $K_s$-saturated graph on $n$ vertices with $\sat(n, K_s)$ edges: the join of a clique with $s-2$ vertices and an independent set with $n-s+2$ vertices.  
\end{theorem}

The \textit{join} $G_1 \ast G_2$ of two graphs $G_1$ and $G_2$ is obtained by taking the disjoint union of $G_1$ and $G_2$ and adding all the edges between them. Erd\H{o}s, Hajnal, and Moon proved Theorem \ref{EHM} by using a clever induction argument. A novel approach to prove this theorem is due to Bollob\'as \cite{B65}, who developed an interesting tool based on systems of intersecting sets. Graph saturation has been studied extensively since Theorem \ref{EHM} appeared half a century ago (see, e.g., \cite{FFS} for a very informative survey). Alon and Shikhelman's generalization of the Tur\'an number motivated Kritschgau, Methuku, Tait, and Timmons \cite{KMTT} to start the systematic study of the function $\sat(n, H, F)$, which denotes the minimum number of copies of $H$ in an $n$-vertex $F$-saturated graph. Here again note that $\sat(n, K_2, F) = \sat(n, F)$. Historically, a natural generalization of counting the number of edges ($K_2$) is to count the number of cliques ($K_r$) of a fixed size, see e.g., \cite{B76}, \cite{E}, and \cite{Z}, where the authors answered the generalized extremal question of finding the maximum number of $K_r$'s in a $K_s$-free graph with fixed number of vertices. Towards generalizing Theorem \ref{EHM} in a similar fashion, Kritschgau, Methuku, Tait, and Timmons proved the following lower and upper bounds, which differ by a factor of about $r-1$, and conjectured that the upper bound (achieved by the same construction given in Theorem \ref{EHM}) is correct.

\begin{theorem} [Kritschgau, Methuku, Tait, and Timmons 2018] \label{tait}
For every $s > r \ge 3$, there exists a constant $n_{r,s}$ such that for all $n \ge n_{r,s}$, 
\begin{align*}
\max \left\{\frac{\binom{s-2}{r-1}}{r-1} \cdot n - 2 \binom{s-2}{r-1}, \frac{\binom{s-2}{r-1} + \binom{s-3}{r-2}}{r} \cdot n\right\} &\le \sat(n, K_r, K_s) 
\\ &\le (n - s + 2) \binom{s-2}{r-1} + \binom{s-2}{r} .
\end{align*}
\end{theorem}

Our first main contribution confirms their conjecture for sufficiently large $n$ by showing that the upper bound is indeed the correct answer. We also show that the natural construction is the unique extremal graph for this generalized saturation problem for large enough $n$. Furthermore, we prove a corresponding stability result for sufficiently large $n$ which shows that even if we allow up to some $cn$ more copies of $K_r$ than $\sat(n, K_r, K_s)$ in an $n$-vertex $K_s$-saturated graph, the extremal graph will still be the same and unique. It is worth noting that there are relatively few stability results in the area of graph saturation, essentially only \cite{AFGS} by Amin, Faudree, Gould, and Sidorowicz, and \cite{BFP} by Bohman, Fonoberova, and Pikhurko. In the notation of joins, the extremal graph in our problem is $K_{s-2} \ast \overline{K}_{n-s+2}$, i.e., the join of a clique with $s-2$ vertices and an independent set with $n-s+2$ vertices. 

\begin{theorem} \label{sat}
For every $s > r \ge 2$, there exists a constant $n_{r,s}$ such that for all $n \ge n_{r,s}$, we have $\sat(n, K_r, K_s) = (n-s+2) \binom{s-2}{r-1} + \binom{s-2}{r}$. Moreover, there exists a constant $c_{r,s} > 0$ such that the only $K_s$-saturated graph with up to $\sat(n, K_r, K_s) + c_{r,s} n$ many copies of $K_r$ is $K_{s-2} \ast \overline{K}_{n-s+2}$.
\end{theorem}

\begin{remark*}
The moreover part of this theorem is tight in the sense that Theorem \ref{sat} fails for $c_{r,s} = \binom{s-3}{r-2}$. To see that consider the graph $G$ on $n$ vertices which is the join of two graphs $G_1$ and $G_2$, where $G_1$ is $K_{s-1}$ minus an edge, and $G_2$ is an independent set on $n-s+1$ vertices. Clearly, $G$ is $K_s$-saturated, with $\left(2 \binom{s-3}{r-2} + \binom{s-3}{r-1}\right) (n-s+1) + 2\binom{s-3}{r-1} + \binom{s-3}{r}$ many copies of $K_r$. 
\end{remark*}

In the process of proving Theorem \ref{sat}, we consider a more general setting and prove and use an intermediate result, which may also be of independent interest. The condition that $G$ is $F$-saturated can be weakened by removing the condition that $G$ is $F$-free (as also studied in \cite{B78} and \cite{T92}). Perhaps counterintuitively, despite the fact that this is a weaker condition, the literature calls $G$ \textit{strongly $F$-saturated} if adding any edge to $G$ creates a new copy of $F$. Following the notation in the literature, we write $\ssat(n, H, F)$ to denote the minimum number of copies of $H$ in an $n$-vertex strongly $F$-saturated graph. It is obvious that $\ssat(n, H, F) \le \sat(n, H, F)$. We have the following asymptotic result for the function $\ssat$ for cliques.    

\begin{theorem} \label{ssat}
For every $s > r \ge 2$, we have $\ssat(n, K_r, K_s) = n \binom{s-2}{r-1} - o(n)$.
\end{theorem}

Kritschgau, Methuku, Tait, and Timmons \cite{KMTT} showed an interesting result, which says that for any natural number $m$, there are graphs $H$ and $F$ such that $\sat(n,H,F) = \Theta(n^m)$. They showed this as an implication of the following bounds that they proved on the minimum number of $C_r$'s in an $n$-vertex $K_s$-saturated graph.

\begin{theorem} [Kritschgau, Methuku, Tait, and Timmons 2018] \label{tait2}
For $s \ge 5$ and $r \le 2s - 4$, $\sat(n, C_r, K_s) = \Theta(n^{\floor{\frac{r}{2}}})$. More precisely,  
\begin{align*}
& \left(1 - o(1)\right) \frac{n^k (s-2)_k}{4 k} \le \sat(n, C_r, K_s) \le \left(1 + o(1)\right) \frac{n^k (s-2)_k}{2 k} &&\text{	if	} 2 \mid r
\\& \left(1 - o(1)\right) \frac{n^k (s-2)_{k+1} (k-2)!}{r (r-3) (r)_k (s-1)} \le \sat(n, C_r, K_s) \le \left(1 + o(1)\right) \frac{n^k (s-2)_{k+1}}{2} &&\text{	if	} 2 \nmid r
\end{align*}
where $k = \floor{\frac{r}{2}}$ and $(m)_k = m(m-1) \cdots (m-k+1)$.
\end{theorem}

Note that the same construction as in Theorem \ref{EHM} proves the upper bound in Theorem \ref{tait2}. We explain the counting for the upper bound in the proof of our Theorem \ref{cycle} in Section 4. We show that for all sufficiently large $n$ the same natural construction is indeed the unique extremal graph for most $r$.

\begin{theorem} \label{cycle}
For every $s \ge 4$ and odd $r$ with $r \ge 7$ or even $r$ with $r \ge 4 \sqrt{s-2}$, there exists a constant $n_{r,s}$ such that for all $n \ge n_{r,s}$, $K_{s-2} \ast \overline{K}_{n-s+2}$ has the minimum number of copies of $C_r$ among $n$-vertex $K_s$-saturated graphs. Moreover, when also $r \le 2s - 4$ this is the unique such graph.
\end{theorem}

We remark here that for any $r,s$ that do not satisfy the assumptions that $s \ge 4$ and $ r \le 2s-4$, we have $\sat(n, C_r, K_s) = 0$, which can be seen from the same extremal graph $K_{s-2} \ast \overline{K}_{n-s+2}$. In Theorem \ref{cycle}, we could write the explicit value of $\sat(n, C_r, K_s)$, which is just the number of cycles of length $r$ in the graph $K_{s-2} \ast \overline{K}_{n-s+2}$. We chose not to do so because the explicit number is not particularly elegant. Also, it turns out that we are able to find the correct asymptotic answers for $r=4$ and $r=5$, which we include in the sections proving Theorem \ref{cycle}.

Next we turn our attention to a long-standing, yet very fundamental conjecture made by Tuza \cite{T86,T88}. In contrast to the Tur\'an number, one of the inherent challenges in studying the saturation number $\sat(n, H)$ for general graphs $H$ is that this function lacks monotonicity properties that one might hope for. For example, Pikhurko \cite{P} showed that there is a pair of connected graphs $F_1 \subset F_2$ on the same vertex set such that $\sat(n, F_1) > \sat(n, F_2)$ for large $n$, violating monotonicity in the second parameter. Regarding non-monotonicity in the first parameter, K\'aszonyi and Tuza \cite{KT} observed that $\sat(2k-1, P_3) = k+1 > k = \sat(2k, P_3)$ where $P_3$ is the path with $3$ edges. Moreover, Pikhurko showed a wide variety of examples of irregular behavior of the saturation function in \cite{P}. All of this non-monotonicity makes proving statements about the saturation function difficult, in particular because inductive arguments generally do not work. However, in order to find some smooth behavior of the saturation function Tuza conjectured the following.

\begin{conjecture} [Tuza 1986] \label{con}
For every graph $F$, the limit $\lim_{n \rightarrow \infty} \frac{\sat(n, F)}{n}$ exists.
\end{conjecture}

Not much progress has been made towards settling the conjecture. The closest positive attempt was made by Truszczy\'nski and Tuza \cite{TT}, who showed that for every graph $F$, if $\lim \inf_{n \rightarrow \infty} \frac{\sat(n, F)}{n} < 1$, then $\lim_{n \rightarrow \infty} \frac{\sat(n, F)}{n}$ exists and is equal to $1 - \frac{1}{p}$ for some positive integer $p$.

Pikhurko considered the saturation number for graph families to make progress in the negative direction of Conjecture \ref{con}. For a family of graphs $\mathcal{F}$, the saturation number $\sat(n, \mathcal{F})$ is defined to be the minimum number of edges in an $n$-vertex $\mathcal{F}$-saturated graph, where a graph $G$ is called $\mathcal{F}$-saturated if $G$ does not contain a copy of any graph in $\mathcal{F}$ and adding any edge to $G$ will create a copy of a graph in $\mathcal{F}$. Pikhurko first showed in \cite{P01} that there exists an infinite family $\mathcal{F}$ of graphs for which $\lim_{n \rightarrow \infty} \frac{\sat(n, \mathcal{F})}{n}$ does not exist, and later in \cite{P} proved the same for a graph family of size only $4$. We make the first progress in 15 years, moving one step closer. 

\begin{theorem} \label{pro}
There exist infinitely many graph families $\mathcal{F}$ of size $3$ such that the ratio $\frac{\sat(n,\mathcal{F})}{n}$ does not converge as $n$ tends to infinity.
\end{theorem}

In the spirit of considering the generalized saturation number, it is natural to ask the more general question of whether $\lim_{n \rightarrow \infty} \frac{\sat(n, K_r, F)}{n}$ exists for every graph $F$. We remark that this problem is interesting since the order of $\sat(n, K_r, \mathcal{F})$ is linear in $n$ for every graph family $\F$, which can easily be shown by considering the same construction used by K\'aszonyi and Tuza in \cite{KT}, who showed the same for $r=2$. We show that our construction of graph families of size $3$ can be extended to this scenario. 

\begin{theorem} \label{construction}
For every $r \ge 3$, there exist infinitely many graph families $\mathcal{F}$ of size $3$ such that the ratio $\frac{\sat(n, K_r, \mathcal{F})}{n}$ does not tend to a limit as $n$ tends to infinity.
\end{theorem}

We next show an example of a graph $F_r$ for which the function $\sat(n, C_r, F_r)$ behaves irregularly. To be precise, we show that for certain $F_r$, the value of the saturation function depends on certain divisibility conditions of $n$, and the sequence $\sat(n, C_r, F_r)$ oscillates.

\begin{theorem} \label{easy_cons}
For every $r \ge 5$, there exists a graph $F_r$ such that $\sat(n, C_r, F_r)$ is zero for infinitely many values of $n$ and also positive infinitely often.
\end{theorem}

The remainder of this paper is organized as follows. In the next section we prove an asymptotically tight lower bound on $\ssat(n, K_r, K_s)$. Then, we use the results and notations of that section to determine $\sat(n, K_r, K_s)$ exactly for sufficiently large $n$ in Section $3$. Next in Section $4$, we prove a few lemmas which will be useful for computing $\sat(n, C_r, K_s)$, i.e., Theorem \ref{cycle}. We handle the cases of even and odd $r$ in Theorem \ref{cycle} separately, and those will be proved in the subsequent two sections. In Section $7$, we construct infinitely many graph families $\mathcal{F}$ of size $3$ for which the ratio $\frac{\sat(n,\mathcal{F})}{n}$ does not converge. We then extend this construction, with the help of Theorem \ref{ssat}, in Section $8$ in order to prove Theorem \ref{construction}. We prove Theorem \ref{easy_cons} in Section $9$. We finish with a few open problems and concluding remarks in Section $10$. 

\section{Asymptotic result for $\boldsymbol{\ssat(n, K_r, K_s)}$}

In this section, we prove Theorem \ref{ssat}. Let $G = (V,E)$ be an $n$-vertex strongly $K_s$-saturated graph such that the number of $K_r$'s in $G$ is $\ssat(n, K_r, K_s)$. Our aim is to find a lower bound on the number of $K_r$ in $G$. Note that if there is an edge $e \in E$ such that $e$ is not in a copy of $K_r$, then $e$ does not contribute to the number of copies of $K_r$. It turns out that a careful analysis of the edges which are in a copy of $K_r$ saves us the required factor of $r-1$ when we compare against the previous best result (Theorem \ref{tait}). So, it is natural to split the edge set $E$ into two parts in the following manner. Let $E_1$ denote the set of edges which are at least in one copy of $K_r$. Let $E_2 = E \setminus E_1$ be the remaining edges in $G$. Now we will prove a simple but powerful lemma which will be useful throughout the current and next sections.

\begin{lemma} \label{trivial}
Every edge of $E_2$ would not be in a copy of $K_s$ even if any non-edge were added to $G$. 
\end{lemma}

\begin{proof}
Fix an arbitrary edge $uv$ of $E_2$ and an arbitrary non-edge $ab$ of $G$. Note that the sets $\{u,v\}$ and $\{a,b\}$ can overlap, but without loss of generality $b \not \in \{u,v\}$. Assume for the sake of contradiction that after adding the missing edge $ab$ we create a copy of $K_s$ containing both $u$ and $v$. Now if we remove the vertex $b$ from the created copy of $K_s$, we will find a copy of $K_{s-1}$ in $G$ which contains both $u$ and $v$. So $uv$ is in a copy of $K_{s-1}$ in $G$, which contradicts the fact that $uv$ is not in a copy of $K_r$, because $r \le s-1$.
\end{proof}

It will be convenient to define a couple of sets which we will use throughout this section and the next section. For $i = 1,2$, let $G_i$ denote the graph on the same vertex set $V$ with the edge set $E_i$. For a graph $H$, it will be convenient to use the notation $d_{H}(v)$ to denote the degree of $v$ in $H$. It will be useful to split the vertices according to their degree in $G_1$, so we define
\begin{equation} 
A = \{v \in V : d_{G_1}(v) \le n^{\frac{1}{3}}\}. \label{large}
\end{equation}

We can observe that $A$ consists of almost all vertices of $G$, i.e., $|A| = n - o(n)$. This is because $|E_1| \le \binom{r}{2} \ssat(n, K_r, K_s) \le \binom{r}{2} \sat(n, K_r, K_s) = O(n)$, where the last equality follows from the upper bound in Theorem \ref{tait}, and so $|V \setminus A| = O(n^{\frac{2}{3}})$. Now our aim is to show that almost every vertex of $A$ is in a copy of $K_{s-1}$ which has only one vertex of $A$. Note that the extremal graph $K_{s-2} \ast \overline{K}_{n-s+2}$ has this property. Formally, we define the following:
\begin{align} 
B = \{v \in A : \exists a_1, \dots, a_{s-2} \in V \setminus A \text{	such that	} v, a_1, \dots, a_{s-2} \text{	induce a copy of 	} K_{s-1}\}. \label{select}
\end{align}

\begin{lemma} \label{main}
Almost all vertices are in $B$, in the sense that $|B| = n - o(n)$.
\end{lemma}

\begin{proof}
Let $R$ denote the set of vertices in $A$ with degree more than $|A| - 2n^{\frac{2}{3}}$ in the induced subgraph of $G_2$ on $A$. Now we claim that $R$ has at most $2rn^{\frac{2}{3}}$ vertices. Assume for the sake of contradiction that $|R| > 2rn^{\frac{2}{3}}$; then with a simple greedy process we will find a copy of $K_r$ in $G_2$. Start with any vertex $v_1 \in R$, and let $R_1 \subseteq R$ denote the set of vertices in $R$ which are neighbors of $v_1$. Clearly, $|R_1| > 2(r-1)n^{\frac{2}{3}}$ because $v_1$ has less than $2n^{\frac{2}{3}}$ non-neighbors in $R$. For $2 \le i \le r$, we continue this process, i.e., at step $i$ we take a vertex $v_i \in R_{i-1}$, and let $R_i \subseteq R_{i-1}$ denote the set of vertices in $R_{i-1}$ which are neighbors of $v_i$. Clearly, $|R_i| > 2(r-i)n^{\frac{2}{3}}$. Now observe that $v_1, v_2, \dots, v_r$ induce a copy of $K_r$ in $G_2$ which is the desired contradiction. So $|R| \le 2rn^{\frac{2}{3}}$. 

Now our aim is to show that $A \setminus R \subseteq B$, which will be sufficient to finish the proof of this lemma. To this end, fix an arbitrary vertex $v \in A \setminus R$. We will first show that there is $w \in A$ such that $vw$ is not an edge of $G$ and there is no $z \in A$ such that $vz$ and $zw$ are both in $E_1$. This is because there are at most $|A| - 2n^{\frac{2}{3}}$ many $E_2$-neighbors of $v$ in $A$ (which follows from the definition of $R$), and in the induced graph of $G_1$ on $A$, there can be at most $n^{\frac{1}{3}} \left(n^{\frac{1}{3}} - 1\right) = n^{\frac{2}{3}} - n^{\frac{1}{3}}$ vertices at distance $2$ from $v$ (which follows from \eqref{large}). So, there are at least $$|A| - 1 - \left(|A| - 2n^{\frac{2}{3}}\right) - n^{\frac{1}{3}} - \left(n^{\frac{2}{3}} - n^{\frac{1}{3}}\right) = n^{\frac{2}{3}} - 1$$ choices for $w$. Fix such a vertex $w$. As $G$ is $K_s$-saturated, if we added the edge $vw$, then we would create a copy of $K_s$. Furthermore, that $K_s$ cannot contain any vertex from $A$ except $v$ and $w$, because if it contained some $z \in A$, then at least one of $vz$ or $zw$ is in $E_2$, contradicting Lemma \ref{trivial}. Hence there is a copy of $K_{s-1}$ induced by $v$ together with $s-2$ vertices from $V \setminus A$, and so $v \in B$. Therefore, $|B| \ge |A| - |R| \ge n - o(n)$.
\end{proof}

\begin{proof}[Proof of Theorem \ref{ssat}]
For an arbitrary vertex $v \in B$, the number of $K_r$'s induced by $v$ together with $r-1$ vertices from $V \setminus A \subseteq V \setminus B$ is at least $\binom{s-2}{r-1}$ from \eqref{select}. So by Lemma \ref{main}, the number of $K_r$'s in $G$ is at least $\binom{s-2}{r-1} |B| = \binom{s-2}{r-1}n - o(n)$. This matches the upper bound from Theorem \ref{tait}, completing the proof of Theorem \ref{ssat}.
\end{proof}

Note that by defining the set $A$ in \eqref{large} optimally, the best lower bound we can achieve with this argument is that $\ssat(n, K_r, K_s) \ge \binom{s-2}{r-1}n - O\left(\sqrt{n}\right)$. Also note that Theorem \ref{ssat} already proves an asymptotically tight lower bound on $\sat(n, K_r, K_s)$, because: 
$$n \binom{s-2}{r-1} - o(n) \le \ssat(n, K_r, K_s) \le \sat(n, K_r, K_s) \le (n-s+2) \binom{s-2}{r-1} + \binom{s-2}{r}.$$

\section{Exact result for $\boldsymbol{\sat(n, K_r, K_s)}$}
In this section, we will find the exact value of $\sat(n, K_r, K_s)$ for all sufficiently large $n$, proving Theorem \ref{sat}. The same argument will also show that the graph $K_{s-2} \ast \overline{K}_{n-s+2}$ is the unique extremal graph. Moreover, we will prove a stability result, i.e., the same graph is also the unique graph among $K_s$-saturated graphs even if we allow up to some $cn$ more copies of $K_r$ than $\sat(n, K_r, K_s)$. We will start with the structural knowledge we developed in the last section and successively deduce more structure to finally reach the exact structure.

Define $c = \frac{1}{4r^2}$ and consider an $n$-vertex $K_s$-saturated graph $G$ with at most $\sat(n, K_r, K_s) + cn$ copies of $K_r$. By defining the sets $A$ and $B$ as in \eqref{large} and \eqref{select} and applying the same arguments we can make the same structural deductions about $G$ as in the last section. In particular, the number of $K_r$'s with one vertex in $B$ and $r-1$ vertices in $V \setminus A$ is at least 
\begin{equation}
n \binom{s-2}{r-1} - o(n). \label{count}
\end{equation}

Next, define 
\begin{equation}
C = \{v \in B : d_{G_1} (v) > s-2\}. \label{problem}
\end{equation}
For $v \in C$, fix $s-2$ neighbors of $v$ in $V \setminus A$ such that those neighbors along with $v$ induce a copy of $K_{s-1}$ in $G$. For each $v \in C$, pick an edge $vw \in E_1$ such that $w$ is not among the $s-2$ fixed neighbors. Note that the same edge $vw$ can be picked at most once more. Each of these particular edges is in $E_1$, hence these edges are contained in some $K_r$, which is not counted in \eqref{count}. After counting for multiplicity, these extra edges will constitute at least an extra $\frac{|C|}{2\binom{r}{2}}$ many copies of $K_r$. Hence, for sufficiently large $n$, $\frac{|C|}{2\binom{r}{2}} \le 2cn$, which implies that $|C| \le \frac{n}{2}$.

So, the set $B \setminus C$ is non-empty for large enough $n$. We will now prove two more structural lemmas.

\begin{lemma} \label{tight}
  Let $v$ be an arbitrary vertex in $B \setminus C$, and suppose $x_1, x_2, \dots, x_{s-2}$ are vertices in $G$ such that $\{v, x_1, \dots, x_{s-2}\}$ induces a copy of $K_{s-1}$. Then for all $u \in V \setminus \{v\}$ such that $uv$ is not an edge, $u$ is adjacent to all of $x_1, \dots, x_{s-2}$. 
\end{lemma}

\begin{proof}
  Since $\{v, x_1, \dots, x_{s-2}\}$ induces $K_{s-1}$ and $s-1 \geq r$, every edge $vx_i$ is in $E_1$. As $v \in B \setminus C$, $v$ has no more $E_1$-edges.
If we add the non-edge $uv$, we must create a copy of $K_s$. If some vertex $w \not \in \{u, v, x_1, \dots, x_{s-2}\}$ participates in the created copy of $K_s$, then we know that $vw$ must be in $E_2$ since $v$ has no more $E_1$-edges, contradicting Lemma \ref{trivial}. So, the only choice for the remaining $s-2$ vertices of the created copy of $K_s$ would be $x_1, \dots, x_{s-2}$.  Thus $u$ must be adjacent to all of $x_1, \dots, x_{s-2}$.
\end{proof}

\begin{lemma} \label{T}
All vertices of $B \setminus C$ have no incident edges from $E_2$.
\end{lemma}

\begin{proof}
Assume for the sake of contradiction that $uv \in E_2$, where $v \in B \setminus C$. Since $G$ is $K_s$-saturated, $u$ is in a copy $S$ of $K_{s-1} \supseteq K_r$. Since $uv \in E_2$, $v \not \in S$ by Lemma \ref{trivial}. Furthermore, $v$ cannot be adjacent to all the vertices in $S$, or else there would be a copy of $K_s$. Similarly, $v$ is in a copy of $K_{s-1}$, and $u$ is not adjacent to the full set of those vertices. Let $a_1, \dots, a_k$, $b_1, \dots, b_k$ and $c_{k+1}, \dots, c_{s-2}$ be distinct vertices such that $\{u, a_1, \dots, a_k, c_{k+1}, \dots, c_{s-2}\}$ and $\{v, b_1, \dots, b_k, c_{k+1}, \dots, c_{s-2}\}$ both induce $K_{s-1}$. The above argument shows that $k \ge 1$. Now we claim that there must be at least two non-edges between $v$ and the set $\{a_1, \dots, a_k\}$, otherwise the neighbors of $v$ in $\{a_1, \dots, a_k\}$ along with $u, v, c_{k+1}, \dots, c_{s-2}$ will induce a clique of order at least $s-1 \ge r$, which contradicts the fact that $uv \in E_2$. Without loss of generality, $v$ is not adjacent to both $a_1$ and $a_2$. Now by applying Lemma \ref{tight} with $v \in B \setminus C$, and $b_1, \dots, b_k, c_{k+1}, \dots, c_{s-2}$ as $x_1, \dots, x_{s-2}$, and $a_1$ as $u$, we see that $a_1$ is adjacent to all of $b_1, \dots, b_k, c_{k+1}, \dots, c_{s-2}$. The same is true of $a_2$. So, $a_1, a_2, b_1, \dots, b_k, c_{k+1}, \dots, c_{s-2}$ induce a copy of $K_s$ in $G$ which is impossible. 
\end{proof}

\begin{proof} [Proof of Theorem \ref{sat}]
Fix any vertex $v \in B \setminus C$. There exists a set $S$ of $s-2$ vertices such that $S \cup \{v\}$ induces a copy of $K_{s-1}$. Since $v \in B \setminus C$, there are no more $E_1$ edges incident to $v$ other than those to $S$. By Lemma \ref{T}, there are no $E_2$ edges either. By Lemma \ref{tight}, every vertex $u \not \in S \cup \{v\}$ must be adjacent to all vertices in $S$. This is already the graph $K_{s-2} \ast \overline{K}_{n-s+2}$, which is $K_s$-saturated, so $G$ is precisely $K_{s-2} \ast \overline{K}_{n-s+2}$.
\end{proof}

\section{Preparation to compute $\boldsymbol{\sat(n, C_r, K_s)}$}
In this section, we state a few lemmas which will be helpful to prove Theorem \ref{cycle} in the subsequent two sections. Our proof is inspired by the proof in \cite{KMTT}. Compared to that paper, we count the number of cycles more carefully to avoid double-counting, which helps us to get the exact answer. We first find asymptotically the number of cycles of length $r$ in the graph $K_{s-2} \ast \overline{K}_{n-s+2}$. Let $k = \floor{\frac{r}{2}}$ throughout the current and next two sections. There are $\binom{n-s+2}{k}$ many independent sets of order $k$ in the subgraph $\overline{K}_{n-s+2}$. If $r$ is even, then for an arbitrary $k$-vertex independent set $A$, the number of copies of $C_r$ containing $A$ is $\frac{(s-2)_k (k-1)!}{2}$, and each copy of $C_r$ is counted exactly once. If $r$ is odd, then for an arbitrary $k$-vertex independent set $A$, the number of copies of $C_r$ containing $A$ is $\frac{(s-2)_{k+1} k!}{2}$, and each copy of $C_r$ is counted exactly once. Furthermore, there is no copy of $C_r$ with more than $k$ vertices in $\overline{K}_{n-s+2}$ because the maximum independent set of $C_r$ has $k$ vertices. Hence, we have the upper bounds of Theorem \ref{tait2}, i.e., 

\begin{align}
& \sat(n, C_r, K_s) \le \frac{(s-2)_k}{2 k} \cdot n^k + O(n^{k-1}) &&\text{	if	} 2 \mid r \label{optimal1}
\\& \sat(n, C_r, K_s) \le \frac{(s-2)_{k+1}}{2} \cdot n^k + O(n^{k-1}) &&\text{	if	} 2 \nmid r \label{optimal2}
\end{align}

We will use the standard notation $\Theta$ in the next few sections. For two functions $f(n)$ and $g(n)$, we call $f(n) = \Theta(g(n))$ if $0 < \lim\inf_{n \rightarrow \infty} \frac{g(n)}{f(n)} \le \lim\sup_{n \rightarrow \infty} \frac{g(n)}{f(n)} < \infty$.

\begin{lemma} \label{ramsey}
For every fixed $l$, there are $\Theta(n^l)$ independent sets of order $l$ in every $n$-vertex $K_s$-free graph.
\end{lemma}

\begin{proof}
Consider an $n$-vertex $K_s$-free graph $G$. It is obvious that the number of independent sets of order $l$ in an $n$-vertex graph is at most $\binom{n}{l} = \Theta(n^l)$. From the most classical result \cite{ES} in Ramsey theory, we know that $R(l,s)$ exists, where $R(l,s)$ denotes the minimum number $N$ such that every graph of order $N$ contains an independent set of order $l$ or a clique of order $s$. So, for each $R(l,s)$-vertex subset $A$ of $G$, the subgraph induced by $A$ must contain an independent set of order $l$ because $A$ does not contain a copy of $K_s$. Now an independent set of order $l$ can be counted at most $\binom{n-l}{R(l,s)-l}$ times. Accounting for multiple-counts, the number of independent sets of order $l$ in $G$ is at least $\frac{\binom{n}{R(l,s)}}{\binom{n-l}{R(l,s)-l}} = \Theta(n^l)$.  
\end{proof}

Next we give an upper bound on the number of edges of any $n$-vertex $K_s$-saturated graph minimizing the number of copies of $C_r$ with $r \le 2s-4$. It is shown in \cite{KMTT} that for every fixed even $r$, there are $o(n^2)$ many edges in an $n$-vertex $K_s$-saturated graph with minimal number of copies of $C_r$. Next we prove the same for all $r \le 2s-4$. We prove a stronger result for odd $r \le 2s-4$, and repeat the proof for even $r$ from \cite{KMTT} for the sake of completion. 

\begin{lemma} \label{upperbound}
For every $n$-vertex $K_s$-saturated graph $G$ minimizing the number of copies of $C_r$, and for any function $f(n)$ such that $f(n) \rightarrow \infty$ as $n \rightarrow \infty$:
\begin{itemize}
\item For odd $r \le 2s-4$, $G$ has $O\left(n f(n)\right)$ many edges.
\item For even $r$, $G$ has $o(n^2)$ many edges.
\end{itemize}
\end{lemma}

\begin{remark*}
In the case of even $r$, if we could prove that $G$ has $o\left(n^{\frac{3}{2}}\right)$ many edges, then we could follow the proof for odd $r$ and would not have the condition $r \ge 4 \sqrt{s-2}$ in Theorem \ref{cycle}. We have briefly mentioned this again in the concluding remarks. 
\end{remark*}

\begin{proof}
\emph{Case 1: $r$ is odd and $r \le 2s-4$.} We can assume that the function $f(n)$ is such that $f(n) = O(\log n)$. Let $G$ be an $n$-vertex $K_s$-saturated graph minimizing the number of copies of $C_r$. For the sake of contradiction, assume that $G$ has more than $n f(n)$ edges. Let $B$ denote the set of all vertices of $G$ with degree more than $f(n)$. A simple counting implies that $\sum_{v \in B} d(v) \ge n f(n)$. To prove Lemma \ref{upperbound}, it is enough to show that for all $v \in B$, there are at least $\Theta\left(n^{k-1} d(v)\right)$ cycles containing $v$. In this case, the total number of cycles will be at least $\Theta\left(\sum_{v \in B} n^{k-1} d(v)\right) \ge \Theta\left(n^k f(n)\right)$, contradicting \eqref{optimal2} for all sufficiently large $n$. To show this, consider a vertex $v \in B$. Consider an arbitrary independent set $I = \{v_1, \cdots, v_{k-1}\}$ of order $k-1$ in $V(G) \setminus \{v\}$. For every $i \in [k-2]$, choose a set $V_{i,i+1}$ of $s-2$ vertices such that adding the edge $v_iv_{i+1}$ would create a copy of $K_s$ on $\{v_i,v_{i+1}\} \cup V_{i,i+1}$. Let $V_1$ denote an empty set if $vv_1$ is an edge, else set it to be a set of $s-2$ vertices such that adding the edge $vv_1$ would create a copy of $K_s$ on $\{v,v_1\} \cup V_1$. Let $U = I \cup V_1 \cup V_{1,2} \cup \cdots \cup V_{k-2,k-1}$. Let $V'$ denote the set of neighbors of $v$ outside of $U$. Note that $|V'| \ge d(v) - ks \ge \frac{1}{2} \cdot d(v)$ for large enough $n$ (remember that $d(v) \ge f(n)$). For each $a \in V'$, we will show the existence of a cycle of length $r$ containing $a$, $v$, and all vertices in $I$, proving that there are at least $\frac{1}{2} \cdot d(v)$ many copies of $C_r$ containing $I$.
\smallskip \\
\emph{Subcase 1: $av_{k-1}$ and $vv_1$ both are edges.} Pick $k$ distinct vertices $u_1, u_1^*, u_1^{**} \in V_{1,2}, u_2 \in V_{2,3}, \cdots, u_{k-2} \in V_{k-2,k-1}$. This is clearly possible because $|V_{i,i+1}| = s-2$ for all $i$ and $k < s-2$. So, $v v_1 u_1 u_1^* u_1^{**} v_2 u_2 v_3 u_3 \cdots v_{k-2} u_{k-2} v_{k-1} a v$ forms a cycle of length $r$ in $G$.
\smallskip \\
\emph{Subcase 2: $av_{k-1}$ is an edge, but $vv_1$ is not an edge.} Pick $k$ distinct vertices $w, w^* \in V_1, u_1 \in V_{1,2}, \cdots, u_{k-2} \in V_{k-2,k-1}$. This is clearly possible because $|V_1| = s-2$, $|V_{i,i+1}| = s-2$ for all $i$, and $k < s-2$. So, $v w w^* v_1 u_1 v_2 u_2 \cdots v_{k-2} u_{k-2} v_{k-1} a v$ forms a cycle of length $r$ in $G$.
\smallskip \\
\emph{Subcase 3: $av_{k-1}$ is not an edge, but $vv_1$ is an edge.} Pick $k-1$ distinct vertices $u_1, u_1^* \in V_{1,2}, u_2 \in V_{2,3}, \cdots, u_{k-2} \in V_{k-2,k-1}$. This is clearly possible because $|V_{i,i+1}| = s-2$ for all $i$ and $k-1 < s-2$. Choose a set $S$ of $s-2$ vertices such that adding the edge $av_{k-1}$ would create a copy of $K_s$ on $\{a,v_{k-1}\} \cup S$. Now as $I$ is an independent set, no vertex from $I \setminus \{v_{k-1}\}$ can be in $S$, so there is a vertex $c \in S$ that is not in the set $I \cup \{v, u_1, u_1^*, \cdots, u_{k-2}\}$. Hence, $v v_1 u_1 u_1^* v_2 u_2 \cdots v_{k-2} u_{k-2} v_{k-1} c a v$ forms a cycle of length $r$ in $G$.
\smallskip \\
\emph{Subcase 4: $av_{k-1}$ and $vv_1$ both are not edges.} Pick $k-1$ distinct vertices $w \in V_1, u_1 \in V_{1,2}, \cdots, u_{k-2} \in V_{k-2,k-1}$. This is clearly possible because $|V_1| = s-2$, $|V_{i,i+1}| = s-2$ for all $i$, and $k-1 < s-2$. Choose a set $S$ of $s-2$ vertices such that adding the edge $av_{k-1}$ would create a copy of $K_s$ on $\{a,v_{k-1}\} \cup S$. Now as $I$ is an independent set, no vertex from $I \setminus \{v_{k-1}\}$ can be in $S$, so there is a vertex $c \in S$ that is not in the set $I \cup \{v, w, u_1, \cdots, u_{k-2}\}$. Hence, $v w v_1 u_1 v_2 u_2 \cdots v_{k-2} u_{k-2} v_{k-1} c a v$ forms a cycle of length $r$ in $G$.
\smallskip \\
From Lemma \ref{ramsey}, we know that there are $\Theta(n^{k-1})$ many independent sets of order $k-1$ in the induced graph $G \setminus \{v\}$ for any vertex $v$, and for each $v \in B$ and such an independent set, we have $\frac{1}{2} \cdot d(v)$ many copies of $C_r$ containing $v$ and the independent set. It is clear that a copy of $C_r$ in $G$ can be counted at most only a constant (depending on $k$) times in this way. So, the number of $C_r$'s in $G$ is at least $\Theta\left(\sum_{v \in B} n^{k-1} d(v)\right) = \Theta\left(n^k f(n)\right)$, contradicting \eqref{optimal2} for all sufficiently large $n$.
\medskip \\
\emph{Case 2: $r$ is even.} 
By Theorem $1^{**}$ in \cite{ES83}, there exists $c, c' > 0$ such that for any graph $G$ with more than $c n^{2 - \frac{2}{r}}$ edges, there exists $$c'n^r \left(\frac{|E(G)|}{n^2}\right)^{\frac{r^2}{4}}$$ copies of $K_{\frac{r}{2},\frac{r}{2}}$. Therefore, if the number of edges of $G$ is $\epsilon n^2$ for some $\epsilon > 0$ and sufficiently large $n$, then there are $\Theta(n^r)$ copies of $C_r$, contradicting \eqref{optimal1}. 
\medskip \\
Since all the cases give contradictions, we are done.
\end{proof}

It is also shown in \cite{KMTT} that for every even $r$, there are $(1 - o(1)) \binom{n}{k}$ many independent sets of order $k$ in an $n$-vertex $K_s$-saturated graph with minimal number of copies of $C_r$, with an application of the Moon-Moser theorem \cite{MM}. Next we prove the same for all $r \le 2s-4$ by using Lemma \ref{upperbound} and the following lemma which is equivalent to the problem appeared in Exercise 40(b) in Chapter 10 of \cite{L}. 

\begin{lemma} \label{mm}
Let $G$ be a graph on $n$ vertices with $\frac{1}{\tau} \binom{n}{2}$ many edges, where $\tau$ is a positive real number. Let $l$ be a positive integer such that $l \le \tau + 1$. Then, the number of independent sets of order $l$ in $G$ is at least $\binom{\tau}{l} \left(\frac{n}{\tau}\right)^l$.
\end{lemma}

\begin{corollary} \label{indep}
For every $n$-vertex $K_s$-saturated graph $G$ minimizing the number of copies of $C_r$ for some $r \le 2s-4$, $G$ has $(1 - o(1)) \binom{n}{k}$ many independent sets of order $k$.
\end{corollary}

\begin{proof}
Consider an $n$-vertex $K_s$-saturated graph $G$ minimizing the number of copies of $C_r$. The number of edges in $G$ is $o(n^2)$ from Lemma \ref{upperbound}, so we can apply Lemma \ref{mm} to conclude that $G$ has $(1 - o(1)) \binom{n}{k}$ many independent sets of order $k$. 
\end{proof}

Notice that the arguments for the even cycles and the odd cycles are bit different in Lemma \ref{upperbound}. It turns out that the proof of Theorem \ref{cycle} for the cases of even and odd $r$ is very different. So, we split the cases in two subsequent sections.

\section{Few copies of $\boldsymbol{C_r}$ in $\boldsymbol{K_s}$-saturated graphs for odd $\boldsymbol{r}$}
Let $G$ be an $n$-vertex $K_s$-saturated graph minimizing the number of copies of $C_r$. Similar to the proof of Theorem \ref{ssat} in Section 2, we define 
\begin{equation}
A = \{ v \in V : d_G(v) \le n^{\frac{1}{3}}\}. \label{redefine}
\end{equation}
We know that $G$ has $O\left(n \log n\right)$ edges from Lemma \ref{upperbound}, so $|A| = n - o(n)$.

Recall that $r \le 2s-4$ and $k = \floor{\frac{r}{2}}$. Consider the collection of independent sets $I$ of order $k$ in $A$ such that for all $v_1, v_2 \in I$, there is no common neighbor of $v_1$ and $v_2$ in $A$. Denote this collection of such independent sets by $\mathcal{I}$. Clearly, there will be $(1 - o(1)) \binom{n}{k}$ independent sets in $\mathcal{I}$. Now consider an arbitrary independent set $I = \{v_1, \cdots, v_k\} \in \mathcal{I}$. For every $i,j \in [k]$, there exists a set $V_{i,j} \subseteq V \setminus A$ of $s-2$ vertices such that adding the edge $v_iv_j$ would create a copy of $K_s$ on $\{v_i,v_j\} \cup V_{i,j}$. Now an easy but cumbersome calculation (similar to the calculation for \eqref{optimal2}) tells us that the number of copies of $C_r$ containing $I$ and $k+1$ vertices from $V \setminus A$ is at least $\frac{(s-2)_{k+1} k!}{2}$, where equality holds if and only if all $V_{i,j}$'s are the same and $v_i, v_j$ do not have any common neighbor in $(V \setminus A) \setminus V_{ij}$. At this point, we can conclude that the upper bound in equation \eqref{optimal2} is asymptotically tight for all odd $r$. Moreover, we can safely say that there are $(1 - o(1)) \binom{n}{k}$ independent sets $I = \{v_1, \cdots, v_k\} \in \mathcal{I}$ for which $V_{i,j}$'s are the same and $v_i, v_j$ do not have any common neighbor in $(V \setminus A) \setminus V_{ij}$, otherwise $G$ will have more copies of $C_r$ than the upper bound in \eqref{optimal2}, which is a contradiction. Let $\mathcal{J}$ denote the collection of independent sets for which the above holds.

\begin{remark*}
Although the statement of Theorem \ref{cycle} is only for odd $r \ge 7$, the above argument actually asymptotically finds the value of $\sat(n, C_r, K_s)$ for $r = 5$ as well.
\end{remark*}

\begin{lemma} \label{maximum_independent}
For odd $r$ with $7 \le r \le 2s-4$, there is an independent set of order $n - o(n)$ in $G$ such that there is a copy $T$ of $K_{s-2}$ in $G$ with the property that every vertex in $T$ is a neighbor of every vertex of the independent set.
\end{lemma}

\begin{proof}
From the fact that $|\mathcal{J}| = \binom{n}{k} - o(n^k)$, we can say that there exist two vertices $u,v \in A$ such that there are $\binom{n}{k-2} - o(n^{k-2})$ independent sets in $\mathcal{J}$ where each of them contains both $u$ and $v$. Let $\mathcal{K}$ denote the collection of independent sets in $\mathcal{J}$ containing both $u$ and $v$. Let $T \subseteq V \setminus A$ be a set of $s-2$ vertices such that adding the edge $uv$ would create a copy of $K_s$ on $\{u,v\} \cup T$. By the definition of $\mathcal{J}$, all the vertices appearing in an independent set in $\mathcal{K}$ should be neighbors of all the vertices in $T$, hence they will form an independent set (because $G$ does not have a copy of $K_s$). For $r \ge 7$, equivalently for $k \ge 3$, it is easy to check that the number of such vertices is $n - o(n)$ (note that this is not true for $k=2$). So, we are done.
\end{proof}

\begin{proof}[Completing the proof.]
Consider the maximum size independent set $I$ in $G$ such that there exists a copy $T$ of $K_{s-2}$ in $G$ such that every vertex in $T$ is a neighbor of every vertex of the independent set. Let $|I| = n - m$. We know that $m = o(n)$ from Lemma \ref{maximum_independent}. Let $S$ denote the set of all vertices outside of $I$ and $T$. For the sake of contradiction, assume that $S$ is non-empty. Now if we let $m' = |S|$, clearly $m' = m - s + 2 = o(n)$. We claim that any $v \in I$ has at least one neighbor in $S$, which will imply that there are at least $n-m$ edges between $I$ and $S$. If there is some $v \in I$ with no neighbor in $S$, then for any $u \in S$, the copy of $K_s$ created by adding the edge $uv$ cannot contain any vertex from $I$ or $S$ except $u$ and $v$, hence $u$ is neighbor of all the vertices of $T$, which in turn tells us that $u$ cannot have any neighbor in $I$, contradicting the maximal choice of $I$. Thus, every vertex in $I$ has at least one neighbor in $S$. 

Let $z$ be the number such that $(z)_k = \sqrt{n(m+k)} \cdot n^{k-1}$. We will show that all the vertices in $S$ can have at most $z$ neighbors in $I$ for sufficiently large $n$. Suppose for contradiction that $v \in S$ has more than $z = o(n)$ neighbors in $I$. We already know that the number of copies of $C_r$ in the induced subgraph on $T \cup I$ is at least $\frac{(s-2)_{k+1}}{2} \cdot (n-m)_k$. Now for any set of $k$ vertices from the neighbors of $v$, there is at least a copy of $C_r$ containing those vertices, together with the vertex $v$ and $k$ vertices from $T$. Clearly, there will be at least $\binom{z}{k}$ such copies of $C_r$. This implies that the number of copies of $C_r$ in $G$ is at least $\frac{(s-2)_{k+1}}{2} \cdot (n-m)_k + \binom{z}{k}$, which contradicts \eqref{optimal2} for all large $n$ because of the following.

\begin{align*}
\frac{(s-2)_{k+1}}{2} &\cdot (n-m)_k + \binom{z}{k} \\
&\ge \frac{(s-2)_{k+1}}{2} \cdot (n-m-k)^k + \frac{(z)_k}{k!} \\
&\ge \frac{(s-2)_{k+1}}{2} \cdot n^k - \frac{(s-2)_{k+1}}{2} \cdot k (m+k) n^{k-1} + \frac{1}{k!} \sqrt{n(m+k)} \cdot n^{k-1} \\
&\ge \frac{(s-2)_{k+1}}{2} \cdot n^k + \frac{1}{2k!} \sqrt{n(m+k)} \cdot n^{k-1} \\
&\ge \frac{(s-2)_{k+1}}{2} \cdot n^k + \Theta\left(n^{k-\frac{1}{2}}\right)
\end{align*}

We now return to our analysis of $S$, which we assumed to be non-empty for the sake of contradiction. If $m \le \log n$, then the number of edges between $S$ and $I$ is at most $(m - s + 2) z < n - m$ for all sufficiently large $n$, which is a contradiction. The only case remaining to handle is when $m > \log n$. For a vertex $v \in S$, pick a set $B$ of $k-1$ vertices from $I$ which are not neighbors of $v$. Fix an order $v_1, v_2, \cdots, v_{k-1}$ of the vertices in $B$, and choose $(s-2)$-element sets $V_i$ such that adding $vv_i$ creates a copy of $K_s$ on $\{v,v_i\} \cup V_i$. Note that as $v$ is not adjacent to all the vertices in $T$ (this follows from the maximality of $I$), $V_i$ cannot be equal to $T$ for any $i$. Consider copies of $C_r$ containing $v, v_1, \cdots, v_{k-1}$ in that order such that the vertices (one or two) between $v$ and $v_1$ are from $V_1$, the vertices (one or two) between $v$ and $v_{k-1}$ are from $V_{k-1}$, and the rest of the vertices are from $T$. Call these cycles \textit{good}. The number of good cycles is at least $\frac{(s-2)_{k+1} k}{2} + 1$. Then there are at least $\left(\frac{(s-2)_{k+1} k}{2} + 1\right) \cdot (k-1)! \ge \frac{(s-2)_{k+1} k!}{2} + 1$ many copies of $C_r$ of good type containing $v$ and all the vertices in $B$. So, if there is no over-counting, then the number of copies of $C_r$ containing $k-1$ vertices from $I$ and one vertex from $S$ is at least $\binom{n - z}{k-1} (m-s+2) \left(\frac{(s-2)_{k+1} k!}{2} + 1\right)$. To show that there is no over-counting, consider a vertex $v \in S$, $k-1$ non-neighbors $v_1, \cdots, v_{k-1}$ of $v$ in $I$ and sets $V_i$ for which adding $vv_i$ creates a copy of $K_s$ on $\{v,v_i\} \cup V_i$. The good cycles containing exactly one vertex from $V_1$ and one from $V_{k-1}$ cannot be counted twice because there is a unique independent set consisting of one vertex in $S$ and $k-1$ in $I$ in this kind of cycle. Now consider a good cycle with two vertices in $V_1$ or $V_{k-1}$. Without loss of generality, the cycle is of the form $vuu'v_1u_1v_2 \cdots u_{k-2}v_{k-1}wv$ with $u,u' \in V_1$. There is again a unique independent set consisting of one vertex in $S$ and $k-1$ in $I$ in this kind of cycle, because $u,v_1, \cdots, v_{k-1}$ cannot be an independent set due to the fact that $u$ and $v_1$ are adjacent. Hence, the total number of copies of $C_r$ in $G$ is at least $\frac{(s-2)_{k+1}}{2} \cdot (n-m)_k + \binom{n - z}{k-1} (m-s+2) \left(\frac{(s-2)_{k+1} k!}{2} + 1 \right)$. Noting that $z = o(n)$ and $m > \log n$, we can see that this contradicts \eqref{optimal2} due to the following:

\begin{align*}
\begin{aligned}[t]
&\frac{(s-2)_{k+1}}{2} \cdot (n-m)_k + \binom{n - z}{k-1} (m-s+2) \left(\frac{(s-2)_{k+1} k!}{2} + 1 \right) \\
&\ge \frac{(s-2)_{k+1}}{2} \cdot n^k - \frac{(s-2)_{k+1}}{2} \cdot k (m+k) n^{k-1} + \frac{(s-2)_{k+1}}{2} \cdot k m (n-z-k)^{k-1} \\
& \; \; \; \; \; \; \; \; \; \; \; \; \; \; \; \; \; \; \; \; \; \; \; \; \; \; \; \; \; \; \; \; \; \; \; \; \; \; \; \; \; \; \; \; \; \; \; \; \; \; \; \; \; \; \; \; \; \; \; \; \; \; \; \; \; \; \; \; \; \; \; \; \; \; \; \; \; \; \; \; \; \; \; \; \; \; \; \; \; \; \; \; \; + \binom{n - z}{k-1} m - \Theta\left(n^{k-1}\right) \\
&\ge \frac{(s-2)_{k+1}}{2} \cdot n^k - \frac{(s-2)_{k+1}}{2} \cdot k m n^{k-1} - \Theta\left(n^{k-1}\right) + \frac{(s-2)_{k+1}}{2} \cdot k m n^{k-1} \\
& \; \; \; \; \; \; \; \; \; \; \; \; \; \; \; \; \; \; \; \; \; \; \; \; \; \; \; \; \; \; \; \; \; \; \; \; \; \; \; \; \; \; \; \; \; \; \; \; \; \; \; \; \; \; \; \; \; \; \; \; \; \; \; \; \; \; \; \; \; \; \; \; \; \; \; \; \; \; \; \; \; \; \; \; \; \; \; \; - \Theta\left(mzn^{k-2}\right) + \Theta\left(mn^{k-1}\right) \\
&\ge \frac{(s-2)_{k+1}}{2} \cdot n^k + \Theta\left(mn^{k-1}\right).
\end{aligned}
\end{align*}
\end{proof}

This shows that $S$ is an empty set and so, $G$ is the union of $I$ (an independent set with maximum size) and $T$ (which is a $K_{s-2}$) where every vertex in $T$ is incident to every vertex in $I$. This finishes the proof of Theorem \ref{cycle} for odd $r$.

\section{Few copies of $\boldsymbol{C_r}$ in $\boldsymbol{K_s}$-saturated graphs for even $\boldsymbol{r}$}
Let $G$ be an $n$-vertex $K_s$-saturated graph minimizing the number of copies of $C_r$. Let $I$ be an arbitrary independent set of order $k = \frac{r}{2}$ in $G$. Our goal is to count the number of copies of $C_r$ in $G$ containing $I$. There are $\frac{(k-1)!}{2}$ circular permutations of the vertices of $I$, accounting for the directional symmetry of a cycle. Fix such an order $v_1, v_2, \cdots, v_k$. For every distinct $i,j \in [k]$, choose a set $V_{i,j}$ of $s-2$ vertices such that adding the edge $v_iv_j$ would create a copy of $K_s$ on $\{v_i,v_j\} \cup V_{i,j}$. For $i \in [k-1]$, we can iteratively pick a common neighbor $u_i$ of $v_i$ and $v_{i+1}$ among the $s-2$ vertices in $V_{i,i+1}$, and finally pick a common neighbor $u_k$ of $v_1$ and $v_k$ from $V_{1,k}$, thus forming a cycle $v_1 u_1 v_2 u_2 \cdots v_k u_k v_1$. Clearly, the number of ways to do this is at least $(s-2)_k$, so there are at least $\frac{(s-2)_k (k-1)!}{2}$ many copies of $C_r$ containing $I$. But there may be over-counting due to the fact that a cycle of length $r$ has two independent set of order $k$. So, to efficiently account for this double-counting, let us define a notion of `essential count'. The idea is to count a copy of $C_r$ containing two independent sets of order $k$ as half, so that the double-counting will make the count exactly one. So, we have two categories of $C_r$ containing $I$, (i) with two independent sets of order $k$, and (ii) with exactly one independent set of order $k$. Now, if there are $x$ copies of $C_r$ containing $I$ of type (i) and $y$ copies of $C_r$ containing $I$ with type (ii), then we say the essential count of the number of copies of $C_r$ containing $I$ is $\frac{x}{2} + y$. 

For a fixed $k$-independent set $I = \{v_1, \cdots, v_k\}$, we now want to find the essential count of the number of copies of $C_r$ containing $I$ in the order $v_1, v_2, \cdots, v_k$. As before, for every distinct $i,j \in [k]$, choose a set $V_{i,j}$ of $s-2$ vertices such that adding the edge $v_iv_j$ would create a copy of $K_s$ on $\{v_i,v_j\} \cup V_{i,j}$. Define the sets $A_j = V_{j,j+1} \setminus \bigcup_{i \neq j} V_{i,i+1}$, where $V_{k,k+1} = V_{k,1}$. Now observe that for all $j$ when we pick a common neighbor $u_j \in V_{j,j+1}$ of $v_j$ and $v_j$ to count the cycle $v_1u_1v_2 \cdots v_ku_kv_1$, the vertices $u_1, u_2, \cdots, u_k$ will form an independent set in $G$ if and only if $u_j \in A_j$ for all $j \in [k]$. So, if $s_j = |A_j|$, then the essential count is at least the following: $$f(s_1, s_2, \cdots, s_k) = \frac{1}{2} \prod_{j = 1}^k s_j + \sum_{\substack{J \subseteq [k] \\ |J| \neq 0}} \prod_{j \notin J} s_j \prod_{j \in J} (s-2-s_j - \iota(J,j)),$$ where $\iota(J,j)$ denotes the number of elements in $J$ smaller than $j$.  

\begin{lemma} \label{calculus}
For any $k \ge 2 \sqrt{s-2}$, the function $f(s_1, s_2, \cdots, s_k)$ attains its minimum uniquely at $(0,0, \cdots, 0)$ over the region $\{0, 1, \cdots, s-2\}^k$.
\end{lemma}

\begin{proof}
For a fixed $j$, if we fix all the variables except $s_j$ and vary $s_j$, $f$ is a linear function with respect to $s_j$. So, the minimum will occur either at $s_j = 0$ or $s_j = s-2$, when other variables are fixed. Hence, applying the same argument for all variables, we can conclude that the minimum can occur only at the vertices of the cube $[0,s-2]^k$. It is easy to check that if we evaluate $f$ at a vertex with at least one co-ordinate $0$ and one co-ordinate $s-2$, then the value will be strictly greater than $f(0,0, \cdots, 0)$. Now, the only thing we need to verify is that $f(s-2,s-2, \cdots, s-2) > f(0,0, \cdots, 0)$, which is equivalent to $\frac{1}{2} (s-2)^k > (s-2)_k$. This holds for $k \ge 2 \sqrt{s-2}$, because:
\begin{align*}
\frac{(s-2)_k}{(s-2)^k} = 1 \left(1 - \frac{1}{s-2}\right) \cdots \left(1 - \frac{k-1}{s-2}\right) 
< e^{-\left(0 + \frac{1}{s-2} + \cdots + \frac{k-1}{s-2}\right)} 
= e^{-\frac{k(k-1)}{2(s-2)}} 
< \frac{1}{2}.
\end{align*}
The function $f$ takes strictly greater values at all vertices in $[0,s-2]^k$ than $f(0,0, \cdots, 0)$, so $f(s_1, s_2, \cdots, s_k)$ is strictly greater than $f(0,0, \cdots, 0)$ for all $(s_1, \cdots, s_k) \neq (0, \cdots, 0)$. As there are finitely many points in $\{0, 1, \cdots, s-2\}^k$, there exists some constant $\epsilon > 0$ (that does not depend on $n$, but may depend on $s$ and $k$) such that $f(s_1, s_2, \cdots, s_k) - f(0, 0, \cdots, 0) \ge \epsilon$ for all $(s_1, \cdots, s_k) \neq (0, \cdots, 0)$.
\end{proof}

\begin{proof}[Completing the proof]
The rest of the proof is similar to the odd $r$ case, and we provide an outline here for completeness. Corollary \ref{indep} and Lemma \ref{calculus} imply that the number of copies of $C_r$ in $G$ is at least $(1-o(1)) \frac{(s-2)_k}{2k} \cdot n^k$, which shows that \eqref{optimal1} is asymptotically tight for $r \ge 4 \sqrt{s-2}$. Now by a similar argument to the odd $r$ case, the number of independent sets $I = \{v_1, \cdots, v_k\}$ of order $k$, for which there is a set $V' \subseteq V(G)$ of size $s-2$ such that for all $i \neq j$, $v$ is a common neighbor of $v_i$ and $v_j$ if and only if $v \in V'$, is $(1-o(1))\binom{n}{k}$. Next we have the following lemma whose proof is the same as Lemma \ref{maximum_independent}.

\begin{lemma} \label{maximum_independent1}
For even $r$ with $4 \sqrt{s-2} \le r \le 2s-4$, there is an independent set $I$ of order $n - o(n)$ in $G$ such that there is a copy $T$ of $K_{s-2}$ in $G$ with the property that every vertex in $T$ is a neighbor of every vertex of $I$.
\end{lemma} 

Define the sets $I$, $T$ and $S$, and the numbers $m$, $m'$ and $z$ as after Lemma \ref{maximum_independent}. Following the proof in the last section, we can show that all the vertices in $S$ can have at most $z$ neighbors in $I$ for sufficiently large $n$. For the sake of contradiction, assume that $S$ is non-empty. As before, the case when $m \le \log n$ leads to a contradiction, and $m > \log n$ remains the only case to resolve. For an arbitrary vertex $v \in S$ and an arbitrary set $B$ of $k-1$ vertices from $I$ that are not neighbors of $I$, consider the good cycles (as defined in the last section) containing $v$ and all the vertices in $B$. Like before, if there is no over-counting, then the number of copies of good $C_r$ containing $k-1$ vertices from $I$ and one vertex from $S$ is at least $\binom{n - z}{k-1} (m-s+2) \left(\frac{(s-2)_k (k-1)!}{2} + 1\right)$. To show that there is no over-counting, it turns out that the situation is simpler in this case compared to the odd $r$ case, which follows from the fact that the good cycles always have a unique independent set consisting of one vertex in $S$ and $k-1$ in $I$. Hence, the total number of copies of $C_r$ in $G$ is at least $\frac{(s-2)_k}{2k} \cdot (n-m)_k + \binom{n - z}{k-1} (m-s+2) \left(\frac{(s-2)_k (k-1)!}{2} + 1 \right) \ge \frac{(s-2)_k}{2k} \cdot n^k + \Theta\left(mn^{k-1}\right)$, contradicting \eqref{optimal1}. So, we have completed the proof of Theorem \ref{cycle}.
\end{proof}

Having completed the proof of Theorem \ref{cycle}, we also solve the problem asymptotically for $r = 4$. 

\begin{proposition}
For every $s \ge 4$, we have the following: 
$$\sat(n, C_4, K_s) = (1 + o(1)) \binom{s-2}{2} \binom{n}{2} = (1 + o(1)) \frac{n^2 (s-2)(s-3)}{4}.$$
\end{proposition}

\begin{proof}
The upper bound follows from \eqref{optimal1}. For the lower bound, consider an $n$-vertex $K_s$-saturated graph $G$ minimizing the number of copies of $C_4$. For a non-edge $uv \in E(G)$, choose a set $T \subseteq V(G)$ of $s-2$ vertices such that adding the edge $uv$ would create a copy of $K_s$ on $\{u,v\} \cup T$. Hence the number of copies of $K_4 \setminus e$ (which is the graph after removing an edge from a complete graph on 4 vertices) containing the non-edge $uv$ is at least $\binom{s-2}{2}$. By Lemma \ref{upperbound}, we can conclude that the number of copies of $K_4 \setminus e$ is at least $(1 - o(1)) \binom{s-2}{2} \binom{n}{2}$ (it is a routine to check that we are not doing any multiple-counting). Hence, $\sat(n, C_4, K_s) = (1 + o(1)) \binom{s-2}{2} \binom{n}{2}$.
\end{proof}

\section{Family of size 3 with non-converging saturation ratio}

In this section, we prove Theorem \ref{pro}. We begin by stating the families of graphs that we will use for the construction. 

\begin{definition} \label{def}
For every positive integer $m \ge 4$, let $\F_m$ be the family of the following three graphs.
\begin{itemize} 
\item Let $B_{m,m}$ be the disjoint union of two copies of $K_m$ plus one edge joining them (often called a ``dumb-bell").
\item Let $V_m$ be a copy of $K_m$ plus two more edges incident to a single vertex of the $K_m$.
\item Let $\Lambda_m$ be a copy of $K_m$ plus a single vertex with exactly two edges incident to the $K_m$.  
\end{itemize} 
\end{definition}

The proof of Theorem \ref{pro} boils down to the fact that the behavior of $\sat(n, \F_m)$ depends on whether or not $n$ is divisible by $m$. The following two lemmas constitute the proof.

\begin{lemma} \label{easy}
For every $n$ divisible by $m$, we have $\sat(n, \F_m) \le \frac{n}{m} \binom{m}{2}$.
\end{lemma}

\begin{proof}
Since $n$ is divisible by $m$, the graph $G$ consisting of the disjoint union of $\frac{n}{m}$ many copies of $K_m$ is clearly $\mathcal{F}_m$-saturated, and the number of edges in $G$ is $\frac{n}{m} \binom{m}{2}$, which proves the result. 
\end{proof}

\begin{lemma} \label{family}
For every $n \ge m \ge 4$ where $n$ is not divisible by $m$, we have $\sat(n, \F_m) \ge \frac{n-m}{m} \left(\binom{m}{2} + 1\right)$.
\end{lemma}

The proof of Lemma \ref{family} will easily follow from the next three lemmas about the structure of $\F_m$-saturated graphs. Let $G$ be an $\F_m$-saturated graph on $n$ vertices. Let $B$ be the set of all vertices of $G$ which are contained in any copy of $K_m$.

\begin{lemma} \label{B}
The subgraph induced by $B$ is only a disjoint union of $K_m$'s.  
\end{lemma}  

\begin{proof}
First, note that no subgraph of $G$ is isomorphic to $F_{m,j}$ for any $j \in \{1, 2, \dots, m-1\}$, where $F_{m,j}$ denotes the union of $2$ copies of $K_m$ overlapping in exactly $j$ common vertices. This is because each $F_{m,j}$ contains a copy of $V_m$ or $\Lambda_m$. As $B$ does not have any copies of $F_{m,j}$ for all $j$, all copies of $K_m$ induced by $B$ are pairwise disjoint. Furthermore, the subgraph of $G$ that $B$ induces is just a disjoint union of $K_m$'s, because any other edge would create a copy of $B_{m,m}$ in $G$. 
\end{proof}

Now let $A$ be the set of all vertices not in $B$. Since the structure in $B$ is so simple, our lower bound will follow by independently lower-bounding the number of edges induced by $A$, and the number of edges between $A$ and $B$. We start with $A$.  

\begin{lemma} \label{A}
The set $A$ has at most $m$ vertices, or $A$ is $K_m$-saturated.
\end{lemma}

\begin{proof}
If $A$ is complete, then the number of vertices in $A$ is at most $m$, or else $G$ contains a copy of $K_{m+1}$, and hence $G$ contains a copy of $\Lambda_m$, which is a contradiction. 

So, suppose $A$ is not complete. We claim that adding any edge to the induced graph on $A$ must create a copy of $K_m$ in $G$. Suppose for the sake of contradiction that there is a non-edge $uv$ with $u,v \in A$ such that adding $uv$ does not create a copy of $K_m$ in $G$. However, it must create a copy of one of the graphs $B_{m,m}$, $V_m$, or $\Lambda_m$, hence one of $u$ or $v$ must be in a copy of a $K_m$ in $G$ (because these three graphs have the property that for all edges $ab$, either $a$ or $b$ is in a copy of $K_m$), which contradicts the definition of $A$. 

Finally, we show that adding any edge to the induced graph on $A$ creates a copy of $K_m$ which entirely lies in $A$. Suppose for the sake of contradiction that there is a non-edge $uv$ in the induced graph on $A$ which, if added, would create a copy of $K_m$ which intersects $B$. Let $w \in B$ be a vertex which lies in a created copy of $K_m$ after adding the edge $uv$. That means that $G$ has the edges $uw$ and $vw$, and the copy of $K_m$ in $B$ containing $w$, together with the edges $uw$ and $vw$, creates a copy of $V_m$. So, this is not possible, and we conclude that the
induced subgraph on $A$ is indeed $K_m$-saturated.
\end{proof}

It only remains to bound the number of edges between $A$ and $B$. We have the following structural lemma.

\begin{lemma} \label{AB}
If $A$ is non-empty, then each copy of $K_m$ in $B$ has at least one edge to $A$. 
\end{lemma} 

\begin{proof}
Assume for the sake of contradiction that there is a copy $U$ of $K_m$ in $B$ which does not have an edge to $A$. Consider arbitrary vertices $u \in U$ and $v \in A$. Then one of the following situations must happen. 
\medskip \\
\emph{Case 1: Adding $uv$ creates a copy of $K_m$.} Then it is easy to check that $U$ has an edge to $A$, which is a contradiction.
\medskip \\
\emph{Case 2: Adding $uv$ creates a copy of $B_{m,m}$ with $uv$ being the middle edge connecting the copies of $K_m$.} This would imply that $v$ is in a copy of $K_m$, which is a contradiction.  
\medskip \\
\emph{Case 3: Adding $uv$ creates a copy of $V_m$ or $\Lambda_m$ with $uv$ being one of the two extra edges outside of the copy of $K_m$.} Then if $uv$ becomes one of the extra edges, the other extra edge should already be there and will connect $U$ and $A$, giving a contradiction.
\medskip \\
Since all the cases give contradictions, we are done. 
\end{proof}

We now combine the previous three lemmas to prove Lemma \ref{family}, which then finishes the proof of Theorem \ref{pro}.

\begin{proof}[Proof of Lemma \ref{family}]
Let $n$ and $m$ satisfy the conditions of Lemma \ref{family}. Clearly $A$ must be non-empty because the number of vertices in $B$ is a multiple of $m$ by Lemma \ref{B}, and so Lemma \ref{AB} implies that there are at least $k$ edges between $B$ and $A$, where $k$ is the number of disjoint copies of $K_m$ in $B$. Now from Lemma \ref{A}, we have two situations. When $A$ has at most $m$ vertices, using Lemmas \ref{B} and \ref{AB}, the number of edges in $G$ is at least $\big\lfloor\frac{n}{m}\big\rfloor \binom{m}{2} + \big\lfloor\frac{n}{m}\big\rfloor \ge \frac{n-m}{m} \left(\binom{m}{2} + 1\right)$. Otherwise, $A$ is $K_m$-saturated, so Theorem \ref{EHM} implies that for all $m \ge 4$ the number of edges in $G$ is at least:

\begin{align*}
& k \binom{m}{2} + k + \left(n - (k+1)m + 2\right) (m-2) \\
& > (km) \frac{m-1}{2} + k + \left(n - (k+1)m\right) \left(\frac{m-1}{2} + \frac{1}{m}\right) \\
& = (km) \frac{m-1}{2} + \left(n - (k+1)m\right) \frac{m-1}{2} + k + \left(n - (k+1)m\right) \frac{1}{m} \\
& = \frac{n-m}{m} \left(\binom{m}{2} + 1\right).
\end{align*} 
This completes the proof. 
\end{proof}

\section{Family of size 3 for generalized saturation ratio}

Inspired by the construction for Theorem \ref{pro}, we extend the construction to prove Theorem \ref{construction}. One of the key challenges is to find an appropriate extension of Theorem \ref{EHM}. Fortunately, our Theorem \ref{ssat} rescues us. We start by stating the families of graphs that we will use for the construction, which are not quite the straightforward generalizations of the families used in Theorem \ref{pro}. For notational brevity, let $r \ge 2$ be a fixed integer for the remainder of this section. 

\begin{definition}
For every positive integer $m \ge 2r^2 + 2r$, let $\F_m$ be the family of the following three graphs.
\begin{itemize} 
\item Let $B_{m,m}$ be the same ``dumb-bell" graph from Definition \ref{def}.
\item Let $V_{m,r}$ be the union of a copy of $K_m$ and a copy of $K_{m-r+1}$ overlapping in exactly one common vertex. 
\item Let $\Lambda_{m,r}$ be a copy of $K_m$ plus a single vertex with exactly $r$ edges incident to the $K_m$.  
\end{itemize} 
\end{definition}

Note that for $r = 2$, we have $\Lambda_{m,2} = \Lambda_m$. However, $V_{m,r}$ is not quite a generalization of $V_m$, and in fact $V_m$ is a subgraph of $V_{m,2}$. We considered $V_m$ instead of $V_{m,2}$ in the case of $r = 2$ to make the analysis simpler and more elegant. So, the above construction actually gives different families of three graphs with non-converging saturation ratio for $r = 2$. 

We proceed to the proof of Theorem \ref{construction}. It turns out that the behavior of $\sat(n, K_r, \F_m)$ is similar to before, i.e., it depends on whether or not $n$ is divisible by $m$. The following two lemmas constitute the proof.

\begin{lemma}
For every $n$ divisible by $m$, we have $\sat(n, K_r, \F_m) \le \frac{n}{m} \binom{m}{r}$.
\end{lemma}

\begin{proof}
The same graph used in the proof of Lemma \ref{easy}, i.e., the disjoint union of $\frac{n}{m}$ many copies of $K_m$, gives us the desired upper bound.  
\end{proof}

\begin{lemma} \label{family'}
For every $n \ge m \ge 2r^2 + 2r$ where $n$ is not divisible by $m$, we have that $\sat(n, K_r, \F_m) \ge \frac{n}{m} \left(\binom{m}{r} + 1\right) - o(n)$.
\end{lemma}

Similarly to Lemma \ref{family}, the proof of Lemma \ref{family'} will follow from the next couple of structural lemmas about $\F_m$-saturated graphs. Let $G$ be an $\F_m$-saturated graph on $n$ vertices. Let $B$ be the set of all vertices of $G$ which are contained in any copy of $K_m$. The subgraph induced by $B$ is only a disjoint union of $K_m$'s, by essentially the same proof as Lemma \ref{B}. Now let $A$ be the set of all vertices not in $B$. Motivated by Lemma \ref{A}, we have the following lemma.

\begin{lemma} \label{A'}
$A$ has at most $m$ vertices, or $A$ is strongly $K_{m-r}$-saturated.
\end{lemma}

\begin{remark*}
This is in contrast to Lemma \ref{A}, which got that the induced graph on $A$ was $K_m$-saturated. Here we only get strongly $K_{m-r}$-saturated (recall that despite its counterintuitive name, strong saturation is a weaker condition), but we can later use our Theorem \ref{ssat} to lower-bound the number of copies of $K_r$ in $A$.
\end{remark*}

\begin{proof}
If $A$ is complete, then the number of vertices in $A$ is at most $m$, or else $G$ contains a copy of $K_{m+1}$, and hence $G$ contains a copy of $\Lambda_{m,r}$, which is a contradiction. So, suppose $A$ is not complete. Fix a non-edge $uv$ in the induced graph on $A$. We consider two cases. 
\medskip
\\ \emph{Case 1: Adding $uv$ would create a copy of $K_m$ in $G$.} We will show that the copy of $K_m$ would lie entirely in $A$, giving the required $K_{m-r}$ in $A$. Indeed, assume for the sake of contradiction that there is a non-edge $uv$ in the induced graph on $A$ which, if added, would create a copy of $K_m$ which intersects $B$. That implies that there is a copy $T$ of $K_{m-1}$ which contains the vertex $u$ and intersects $B$. Clearly $T$ can intersect only a single copy $U$ of $K_m$ in $B$, because the induced graph on $B$ is just a disjoint union of $K_m$'s. Now, if $|T \cap U| \ge r$, then $T \cup U$ contains a copy of $\Lambda_{m,r}$, which is a contradiction. Otherwise, $|T \cap U| < r$, and so $T \cup U$ contains a copy of $V_{m,r}$, which is also a contradiction.
\medskip
\\ \emph{Case 2: Adding $uv$ would not create a copy of $K_m$ in $G$.} If adding $uv$ creates a copy of $B_{m,m}$ or $\Lambda_{m,r}$ in $G$, then one of $u$ or $v$ must be in a copy of a $K_m$ in $G$, which contradicts the definition of $A$. Alternatively, if adding $uv$ creates a copy of $V_{m,r}$ in $G$, then that copy of $V_{m,r}$ would contain a copy of $K_m$ in $B$, together with $m-r$ vertices in $A$. These $m-r$ vertices would clearly induce a copy of $K_{m-r}$ after adding $uv$. Hence we are done.
\end{proof}

Next, following the proof of Lemma \ref{family}, we bound the number of $K_r$'s that intersect both $A$ and $B$.

\begin{lemma} \label{AB'}
Suppose $m \ge 2r + 1$. If $A$ is non-empty, then for each copy $U$ of $K_m$ in $B$, there is at least one copy of $K_r$ intersecting both $U$ and $A$.
\end{lemma}

\begin{proof}
Assume for the sake of contradiction that there is a copy $U$ of $K_m$ in $B$ for which there is no copy of $K_r$ intersecting both $U$ and $A$. Consider arbitrary non-adjacent vertices $u \in U$ and $v \in A$. One of the following situations must happen. 
\medskip \\
\emph{Case 1: Adding $uv$ creates a copy $T$ of $K_{m-r}$.} Then it is easy to check that there is a copy of $K_{m-r-1}$ (and hence a copy of $K_r$ if $m \ge 2r + 1$) intersecting both $U$ and $A$, which is a contradiction.
\medskip \\
\emph{Case 2: Adding $uv$ creates a copy of $B_{m,m}$ with $uv$ being the middle edge connecting the copies of $K_m$.} This case is exactly the same as before, i.e., $v$ is in a copy of $K_m$, which is a contradiction.  
\medskip \\
\emph{Case 3: Adding $uv$ creates a copy of $\Lambda_{m,r}$ with $uv$ being one of the $r$ extra edges outside of the copy of $K_m$.} Then if $uv$ becomes one of the extra $r$ edges, the $r-1$ endpoints in $U$ of the remaining $r-1$ extra edges, together with the vertex $v$, induce a copy of $K_r$, giving a contradiction.
\medskip \\
Since all the cases give contradictions, we are done. 
\end{proof}

\begin{proof}[Proof of Lemma \ref{family'}]
  Let $n$ and $m$ satisfy the conditions of Lemma \ref{family'}. Clearly $A$ must be non-empty because the number of vertices in $B$ is a multiple of $m$, and so Lemma \ref{AB'} implies that there are at least $k$ copies of $K_r$ intersecting both $B$ and $A$, where $k$ is the number of disjoint copies of $K_m$ in $B$. Now from Lemma \ref{A'}, we have two situations. When $A$ has at most $m$ vertices, by Lemma \ref{AB'}, the number of copies of $K_r$ in $G$ is at least $\big\lfloor\frac{n}{m}\big\rfloor \binom{m}{r} + \big\lfloor\frac{n}{m}\big\rfloor \ge \frac{n-m}{m} \left(\binom{m}{r} + 1\right)$. Otherwise, $A$ is strongly $K_{m-r}$-saturated, so Theorem \ref{ssat} implies that $A$ induces at least $\binom{m-r-2}{r-1} (n-km) - o(n)$ many copies of $K_r$, and so for all $m \ge 2r^2 + 2r$ the number of copies of $K_r$ in $G$ is at least: 
\begin{align}
k \binom{m}{r} + k + & \binom{m-r-2}{r-1} (n-km) - o(n). \label{binomial}
\end{align}
To get the required lower bound, we next prove the simple claim that $\binom{m-r-2}{r-1} \ge \frac{1}{m} \left(\binom{m}{r} + 1\right)$ for all $m \ge 2r^2 + 2r$ and $r \ge 2$. The most convenient way to do this is to show that $m \binom{m-r-2}{r-1} > \binom{m}{r}$, since both sides of this last inequality are integers. Indeed, let $m$ and $r$ satisfy the conditions we just mentioned. Then,
$$\frac{m-1}{m-r-2} \le \frac{m-2}{m-r-3} \le \cdots \le \frac{m-r+1}{m-2r} \le \frac{2r^2 + r + 1}{2r^2}.$$
Hence, $$\frac{\binom{m}{r}}{m \binom{m-r-2}{r-1}} \le \frac{1}{r} \left(1 + \frac{r + 1}{2r^2}\right)^{r-1} \le \frac{1}{r} \cdot e^{\frac{(r+1)(r-1)}{2r^2}} \le \frac{1}{r} \cdot \sqrt{e} < 1,$$ 
which establishes the claim that $\binom{m-r-2}{r-1} \ge \frac{1}{m} \left(\binom{m}{r} + 1\right)$. Using this, we get that \eqref{binomial} is at least $k \binom{m}{r} + k + \frac{1}{m} \left(\binom{m}{r} + 1\right) (n-km) - o(n) \ge \frac{n}{m} \left(\binom{m}{r} + 1\right) - o(n)$, completing the proof.
\end{proof}

\section{Irregular behavior of $\boldsymbol{\sat(n, C_r, F_r)}$}
In this section, we prove Theorem \ref{easy_cons}. In particular, we will prove that for every $r \ge 4$, $\lim \inf_{n \rightarrow \infty} \sat(n, C_{r+1}, B_{r,r}) = 0$, and $\lim \sup_{n \rightarrow \infty} \sat(n, C_{r+1}, B_{r,r}) > 0$, where $B_{r,r}$ is the same ``dumb-bell'' graph from Definition \ref{def}. We remark here that this statement is false for $r = 2$ and $r = 3$, which we show in Proposition \ref{excess} at the end of this section. The following two lemmas constitute the entire proof of Theorem \ref{easy_cons}.

\begin{lemma} \label{si}
For every $n$ divisible by $r$, we have $\sat(n, C_{r+1}, B_{r,r}) = 0$.
\end{lemma}

\begin{proof}
The same graph used in the proof of Lemma \ref{easy}, i.e., the disjoint union of $\frac{n}{r}$ many copies of $K_r$, is $B_{r,r}$-saturated but has no copies of $C_{r+1}$, proving the lemma.  
\end{proof}

\begin{lemma} \label{cycle_sat}
For every $n \ge 2r$, and $r \ge 4$ such that $n$ is not divisible by $r$, we have $\sat(n, C_{r+1}, B_{r,r}) \ge 1$.
\end{lemma}

\begin{proof}
Let $G$ be a $B_{r,r}$-saturated graph on $n$ vertices. We show that there is a cycle of length $r+1$ in $G$ if the conditions of Lemma \ref{cycle_sat} are met. We divide the proof in two cases.
\medskip \\
\emph{Case 1: There is a copy of $F_{r,j}$ in $G$ for some $j \in \{1, 2, \cdots, r-1\}$, where $F_{r,j}$ denotes the union of $2$ copies of $K_r$ overlapping in exactly $j$ common vertices.} It is easy to check that $F_{r,j}$ contains a copy of $C_{r+1}$ for every $j \ge 2$. Hence, if $G$ contains a copy of $F_{r,j}$ for some $j \ge 2$, then there is already a cycle of length $r+1$ in $G$. So, we can assume that $G$ contains a copy of $F_{r,1}$. Assume that $w, u_1, u_2, \cdots, u_{r-1}, v_1, v_2, \cdots, v_{r-1}$ are distinct vertices such that $\{w, u_1, \cdots, u_{r-1}\}$ and $\{w, v_1, \cdots, v_{r-1}\}$ both induce $K_r$. Note that if there is an edge $u_iv_j$ for some $i,j$, then it is easy to find a copy of $C_{r+1}$ using the edge $u_iv_j$. For example, if $u_1v_1$ is an edge, then $w v_1 u_1 u_2 \cdots u_{r-1} w$ forms a $C_{r+1}$. So, we can assume that there is no edge $u_iv_j$ for any $i,j$. Now one of the following situations must happen.
\medskip \\
\emph{Subcase 1: Adding $u_1v_1$ creates a copy of $K_r$.} So, $u_1$ and $v_1$ must have at least $r-2$ common neighbors. If $r \ge 4$, then among $r-2 \ge 2$ common neighbors of $u_1$ and $v_1$, we can pick a vertex $x$ which is distinct from $w$. Now it is easy to check that $w v_1 x u_1 u_2 \cdots u_{r-2} w$ forms a cycle of length $r+1$. 
\medskip \\
\emph{Subcase 2: Adding $u_1v_1$ creates a copy of $B_{r,r}$ with $u_1v_1$ being the middle edge connecting the copies of $K_m$.} Hence, there is either a copy of $K_r$ containing $u_1$ and not containing any vertex in $\{w, v_1, \cdots, v_{r-1}\}$, or a copy of $K_r$ containing $v_1$ and not containing any vertex in $\{w, u_1, \cdots, u_{r-1}\}$. Due to symmetry, it is enough to check the first situation. If there is a copy of $K_r$ containing $u_1$ and not containing any vertex in $\{w, v_1, \cdots, v_{r-1}\}$, then that copy of $K_r$ along with the copy of $K_r$ induced by $\{w, v_1, \cdots, v_{r-1}\}$ and the edge $u_1w$ forms a copy of $B_{r,r}$, which is a contradiction.
\medskip \\
\emph{Case 2: There is no copy of $F_{r,j}$ in $G$ for any $j \in \{1, 2, \cdots, r-1\}$.} Let $B$ be the set of all vertices of $G$ which are contained in any copy of $K_r$. Firstly note that $B$ cannot be empty, because there are two disjoint copies of $K_r$ in the graph $B_{r,r}$ and it is not possible to create two disjoint copies of $K_r$ by adding one edge to $G$. The subgraph induced by $B$ is only a disjoint union of $K_r$'s, by the same proof as Lemma \ref{B}. Now let $A$ be the set of all vertices not in $B$. Clearly $A$ must be non-empty, because the number of vertices in $B$ is a multiple of $r$, and $r$ does not divide $n$. Fix a copy $T = \{v_1, v_2, \cdots, v_r\}$ of $K_r$ in $B$. Now one of the following situations must happen.
\medskip \\
\emph{Subcase 1: There is at most one vertex in $T$ which has edges to $A$.} If there is no edge between $T$ and $A$, then some easy case-checking (similar to before) implies that adding an edge $v_1 a$ (where $a \in A$) would not create a copy of $B_{r,r}$. Now assume that there is exactly one vertex (without loss of generality $v_1$) in $T$ which has edges to $A$. Again some easy case-checking will tell us that adding $v_2 a$ for any $a \in A$ would not create a copy of $K_r$ (because $r \ge 4$), so, the only way to create a copy of $B_{r,r}$ would be to become the middle edge connecting the copies of $K_r$ in $B_{r,r}$, but that would contradict the fact that $a \in A$ (remember that no vertices in $A$ are in a copy of $K_r$). These are all contradictions.
\medskip \\
\emph{Subcase 2: There are at least two vertices in $T$ which have edges to $A$.} If there is $a \in A$ such that $a$ is adjacent to at least two vertices in $T$, then one can find a cycle of length $r+1$ (for example, without loss of generality $a$ is adjacent to both $v_1$ and $v_2$, so, $v_1 a v_2 v_3 \cdots v_r v_1$ forms a copy of $C_{r+1}$). So, we can assume that for all $a \in A$, the vertex $a$ is adjacent to at most one vertex in $T$. Without loss of generality, $v_1, v_2 \in T$ have edges to $A$. Let $v_1 a_1$ and $v_2 a_2$ be edges for some $a_1, a_2 \in A$. If $a_1 a_2$ is an edge, then $v_1 a_1 a_2 v_2 v_3 \cdots v_{r-1} v_1$ forms a cycle of length $r+1$. Now if $a_1 a_2$ is a non-edge, then adding $a_1 a_2$ would create a copy of $K_r$ (because it must create a copy of $B_{r,r}$, but $a_1 a_2$ could not become the middle edge in the created copy of $B_{r,r}$ due to the fact that $a_1$ is not in a copy of $K_r$ in $G$). Note that $a_1$ and $a_2$ cannot have a common neighbor in $T$, so they must have a common neighbor $x \not \in T$, which implies that $v_1 a_1 x a_2 v_2 \cdots v_{r-2}$ forms a cycle of length $r+1$ in $G$.
\medskip \\
Since all the cases either find a cycle of length $r+1$ or give contradictions, we are done.
\end{proof}

\begin{proposition} \label{excess}
For $r = 2$ and $r = 3$, we have $\sat(n, C_{r+1}, B_{r,r}) = 0$ for all $n \ge 2r$.
\end{proposition}

\begin{proof}
In Lemma \ref{si}, we have already seen that $\sat(n, C_{r+1}, B_{r,r}) = 0$ when $n$ is divisible by $r$. So, we have to prove Proposition \ref{excess} when $n$ is not divisible by $r$.

The graph $B_{2,2}$ is the path with 3 edges, i.e., $P_3$. A graph which is a disjoint union of $P_2$'s and $P_1$'s is always $P_3$-saturated. If $n \ge 4$ is odd, then the graph consisting of the disjoint union of a copy of $P_2$ and $\frac{n-3}{2}$ many copies of $P_1$ is a $P_3$-saturated graph with no copy of $C_3$. So, we have $\sat(n, C_3, B_{2,2}) = 0$ for all $n \ge 4$.

For $r = 3$, we split into two cases depending on the value of $n \pmod 3$. If $n$ is of the form $3m + 1$ for some integer $m$, then the graph consisting of $m$ disjoint copies of $K_3$ together with $m$ edges connecting an extra vertex to each copies of $K_3$, is a $B_{3,3}$-saturated graph without any copy of $C_4$. So, for $n \equiv 1 \pmod 3$, we have $\sat(n, C_4, B_{3,3}) = 0$. 

Now when $n$ is of the form $3m + 2$, we have a similar construction. Consider a graph $G$ on the vertex set $\{a,b\} \cup \{x_1, x_2, \cdots, x_m\} \cup \{y_1, \cdots, y_m\} \cup \{z_1, \cdots, z_m\}$, and the edge set $\{x_jy_j : j \in [m]\} \cup \{y_jz_j : j \in [m]\} \cup \{z_jx_j : j \in [m]\} \cup \{ax_j : j \in [m]\} \cup \{by_j : 2 \le j \le m\} \cup \{bx_1\}$. It is easy to verify that $G$ is $B_{3,3}$-saturated, and does not have a copy of $C_4$. Hence, we have $\sat(n, C_4, B_{3,3}) = 0$ for all $n \ge 6$.
\end{proof}

\section{Concluding remarks}
We end with some open problems. We determined the exact value of $\sat(n, K_r, K_s)$ for all sufficiently large $n$, but our arguments do not extend to find the value for small $n$. So, the following question still remains open.

\begin{problem}
For $s > r \ge 3$, determine the exact value of $\sat(n, K_r, K_s)$ for all $n$.
\end{problem}

We have already made a remark on the maximum constant $c_r$ we can write in the stability result in Theorem \ref{sat}. It would be interesting to determine that maximum constant.

\begin{problem}
For $s > r \ge 3$, what is the second smallest number of copies of $K_r$ in an $n$-vertex $K_s$-saturated graph?
\end{problem}

It might be interesting to consider a more general problem of finding the spectrum (set of possible values) of the number of copies of $K_r$ in a $K_s$-saturated graph. The $r = 2$ case, i.e., the edge spectrum of $K_s$-saturated graphs, was completely solved in \cite{AFGS} and \cite{BCFF}.

\begin{problem}
For $s > r \ge 3$, what are the possible numbers of copies of $K_r$ in an $n$-vertex $K_s$-saturated graph?
\end{problem}

We could not extend our method to find the exact value of $\sat(n, C_r, K_s)$ for the situations when $r = 5$, and when $r$ is an even number with $r = O(\sqrt{s})$. So, it will be interesting to find the values of $\sat(n, C_r, K_s)$ for all $r$. We conjecture that the extremal graph $K_{s-2} \ast \overline{K}_{n-s+2}$ is the unique graph minimizing the number of cycles of length $r$ among all $n$-vertex $K_s$-saturated graphs. It is worth mentioning that in Lemma \ref{upperbound}, if we can prove that any $n$-vertex $K_s$-saturated graph with the minimal number of copies of $C_r$ has $o(n^{\frac{3}{2}})$ edges for even $r$, then the proof of Theorem \ref{cycle} for odd $r$ in Section 5 can be adapted for even $r$ as well, and it will prove our conjecture for even $r \ge 6$. 

\begin{problem}
For every $s \ge 4$ and $r \le 2s - 4$, compute the exact value of $\sat(n, C_r, K_s)$.
\end{problem}

Theorem \ref{sat} and Theorem \ref{cycle} motivate us to ask the following general question.

\begin{problem}
Is there a graph $F$, for which $K_{s-2} \ast \overline{K}_{n-s+2}$ does not (uniquely) minimize the number of copies of $F$ among $n$-vertex $K_s$-saturated graphs for all sufficiently large $n$?
\end{problem}

As we mentioned earlier, Conjecture \ref{con} is still wide open and likely needs new ideas to settle it. It would be interesting to figure out if the size of the family in Theorem \ref{pro} can be further reduced to $2$. Finally, as we briefly discussed before stating Theorem \ref{construction}, it would be interesting to consider Conjecture \ref{con} for the generalized saturation problem. 

\begin{problem}
For $r \ge 2$, does the limit $\lim_{n \rightarrow \infty} \frac{\sat(n, K_r, F)}{n}$ exist for every graph $F$?
\end{problem}

\section*{Acknowledgements}
The authors are grateful to the anonymous referees for their suggestions and comments to improve the exposition of this paper. In particular, we are thankful to them for pointing out a technical issue in Lemma \ref{upperbound} in an earlier version of this paper.

\end{document}